\documentclass[11pt]{article}

\usepackage[margin=1in]{geometry}
\usepackage{amsmath,amssymb,amsthm,mathtools}
\usepackage{enumitem}
\usepackage[hidelinks]{hyperref}
\usepackage[nameinlink,capitalize]{cleveref}

\newtheorem{theorem}{Theorem}[section]
\newtheorem{lemma}[theorem]{Lemma}
\newtheorem{corollary}[theorem]{Corollary}
\newtheorem{proposition}[theorem]{Proposition}
\newtheorem{hypothesis}[theorem]{Hypothesis}
\newtheorem{conjecture}[theorem]{Conjecture}
\theoremstyle{remark}
\newtheorem{remark}[theorem]{Remark}
\newtheorem{definition}[theorem]{Definition}
\newtheorem{problem}[theorem]{Problem}
\theoremstyle{plain}

\newcommand{\F}{\mathbb{F}}
\newcommand{\Fp}{\mathbb{F}_\ell}
\newcommand{\Fps}{\mathbb{F}_\ell^{\!*}}

\newcommand{\Z}{\mathbb{Z}}

\newcommand{\Q}{\mathbb{Q}}
\newcommand{\U}{\mathcal U}
\newcommand{\C}{\mathbb C}

\crefname{lemma}{Lemma}{Lemmas}
\Crefname{lemma}{Lemma}{Lemmas}
\crefname{proposition}{Proposition}{Propositions}
\Crefname{proposition}{Proposition}{Propositions}
\crefname{corollary}{Corollary}{Corollaries}
\Crefname{corollary}{Corollary}{Corollaries}
\crefname{theorem}{Theorem}{Theorems}
\Crefname{theorem}{Theorem}{Theorems}
\crefname{definition}{Definition}{Definitions}
\Crefname{definition}{Definition}{Definitions}
\crefname{problem}{Problem}{Problems}
\Crefname{problem}{Problem}{Problems}
\crefname{hypothesis}{Hypothesis}{Hypotheses}
\Crefname{hypothesis}{Hypothesis}{Hypotheses}
\crefname{conjecture}{Conjecture}{Conjectures}
\Crefname{conjecture}{Conjecture}{Conjectures}
\crefname{remark}{Remark}{Remarks}
\Crefname{remark}{Remark}{Remarks}

\title{\huge\bfseries Improved deterministic integer factorization by quadratic-image sieving}
\author{\Large\bfseries  Pantelimon St\u{a}nic\u{a}$^1$
\vspace{.2cm}\\
$^1$Applied Mathematics Department, Naval Postgraduate School\\
Monterey, CA 93943, USA, \texttt{pstanica@nps.edu}\vspace{.15cm}
}
\date{\today}

\begin{document}
\maketitle

\begin{abstract}
Let $N=pq$ be a semiprime and let $\ell\nmid N$ be an odd prime.  For $a\ne0$ modulo $\ell$, put $H_a(N;\ell)=\{ax+Nx^{-1}:x\in\F_\ell^*\}$.  For a primitive root $\gamma$, the size of $H_\gamma(N;\ell)\cup H_{\gamma^2}(N;\ell)$ is determined by the Frobenius trace of $Y^2=(X^2-1)(X^2-\gamma)$.  We study the resulting trace-selection problem through finite certificates, a congruence modulo $8$, and a GRH-conditional positive-density family from the CM specialization $\gamma=2$.

Independently of this conditional result, we obtain a deterministic factorization algorithm from the residue-known coefficient lattice of Harvey and Hittmeir (\emph{Math. Comp.}, 2022).  The actual babystep index is a nondegenerate quadratic image modulo each selected prime divisor of the auxiliary modulus.  After a computable translation, every such image avoids a fixed set of $(\ell-\chi(-1))/4$ residues.  This forbidden set is optimal among translated sets that must handle both square classes of the leading coefficient.  Gauss reduction, a nearby lattice-vector selection, arbitrary-point multipoint evaluation, and the unrestricted large-order algorithm of Harvey and Hittmeir (arXiv, 2026) give the bit-complexity bound
$N^{1/5}\exp(-(\log(4/3)/50+o(1))\log N/\log\log N)$.
Its secondary factor is asymptotically smaller than every fixed negative power of $\log N$, improving the general bound of Harvey and Hittmeir (\emph{Math. Comp.}, 2022) and the balanced-semiprime bound of Gao et al.\ (arXiv, 2025).
\end{abstract}



\section{Introduction}
\label{sec:intro}

Deterministic integer factorization has a long line of improvements in the exponent of $N$.  Harvey proved an exponent-one-fifth deterministic factorization algorithm~\cite{Harvey2021}, and Harvey and Hittmeir later obtained a logarithmic refinement~\cite{HarveyHittmeir2022}.  These algorithms build on Hittmeir's baby-step--giant-step strategy and on Lehman-type coefficient searches combined with local sieving~\cite{Hittmeir2018,Harvey2021,HarveyHittmeir2022,Hittmeir2023}.  The local input is that if $N=pq$, if $\ell\nmid N$, and if $a\in\Fps$, then $ap+q$ modulo $\ell$ belongs to a subset of $\Fp$ occupying roughly half of the residues.  Taking a CRT product over auxiliary primes converts these local restrictions into a possible reduction of a global exponent search.

Harvey and Hittmeir's general deterministic bound is $O(N^{1/5}\log^{16/5}N\,\allowbreak(\log\log N)^{-3/5})$~\cite{HarveyHittmeir2022}.  Gao, Feng, Hu, and Pan recently reduced the logarithmic numerator to $\log^{13/5}N$ for balanced semiprimes~\cite{GaoFengHuPan2025}.  The bound proved here keeps the exponent $1/5$ and replaces these fixed powers of $\log N$ by the factor $\exp(-(\log(4/3)/50+o(1))\log N/\log\log N)$, which is asymptotically smaller than every fixed negative power of $\log N$.

Throughout the paper, $\ell>2$ is an odd prime, $\Fp=\F_\ell$, $\Fps=\Fp\setminus\{0\}$, and $\chi$ is the Legendre symbol on $\Fp$, extended by $\chi(0)=0$.  Thus $\chi(u)=1$ for nonzero squares, $\chi(u)=-1$ for nonsquares, and $\chi(0)=0$.  A primitive root modulo $\ell$ is a generator of the cyclic group $\Fps$.  A semiprime is an integer $N=pq$ with $p$ and $q$ prime.  Locally modulo $\ell$, only the residue class of $N$ is used, and we assume $N\in\Fps$.  For $a\in\Fps$ define $H_a(N;\ell)=\{ax+Nx^{-1}:x\in\Fps\}\subseteq\Fp$.  If $N=pq$ and $\ell\nmid pq$, then taking $x=p$ gives $ap+q\in H_a(N;\ell)$.

We use the following notation and standard facts throughout.  The symbol $\varphi$ denotes Euler's totient function and $\mu$ denotes the M\"obius function.  The Chinese remainder theorem is abbreviated as CRT.  For a statement $P$, the indicator $\mathbf 1[P]$ equals $1$ if $P$ is true and $0$ otherwise.  The logarithmic integral is $\operatorname{Li}(X)=\int_2^X dt/\log t$ for $X>2$.  We use $O(\cdot)$, $\ll$, $\gg$, and $\asymp$ in their usual analytic-number-theory sense, with implied constants independent of the variable tending to infinity, and $f(X)\sim g(X)$ means $f(X)/g(X)\to1$.  For $m\ge1$, $(\Z/m\Z)^*$ denotes the group of residue classes invertible modulo $m$.  If $G$ is a finite abelian group, then $\widehat G$ denotes its group of complex-valued characters.  For $G=\F_\ell^*$, the symbol $\chi_0$ denotes the principal character, and $\operatorname{ord}(\psi)$ denotes the order of a character $\psi\in\widehat{\F_\ell^*}$.  We write $\mathbf 1_{\mathrm{prim}}$ for the indicator of primitive roots in $\F_\ell^*$.  Finally, $\zeta_n$ denotes a fixed primitive $n$th root of unity.

We will use two elementary analytic tools in explicit finite-sum form.  First, if $G$ is a finite abelian group, then character orthogonality says that
\[
  \frac{1}{|G|}\sum_{\psi\in\widehat G}\psi(x)\overline{\psi(y)}
  =
  \begin{cases}
    1,& x=y,\\
    0,& x\ne y,
  \end{cases}
  \qquad x,y\in G.
\]
For a function $f:G\to\C$, define its Fourier transform by $\widehat f(\psi)=\sum_{x\in G}f(x)\overline{\psi(x)}$.  The finite-group Plancherel identity is
\[
  \sum_{\psi\in\widehat G}|\widehat f(\psi)|^2
  = |G|\sum_{x\in G}|f(x)|^2.
\]
For a subset $A\subseteq G$, we use only the specialization $\sum_{\psi\in\widehat G}|\sum_{x\in A}\psi(x)|^2=|G|\,|A|$.  These identities are standard in Fourier analysis on finite abelian groups, as described in Terras~\cite[Part~I]{Terras1999}.  Second, for finite complex sequences $(u_i)$ and $(v_i)$, the Cauchy--Schwarz inequality is
\[
  \left|\sum_i u_i v_i\right|
  \le
  \left(\sum_i |u_i|^2\right)^{1/2}
  \left(\sum_i |v_i|^2\right)^{1/2}.
\]
We use this form in the character-sum estimate below, as stated in Steele~\cite[Chapter~1]{Steele2004}.

We also recall the point-counting estimate used below.  The genus-one case of the Hasse--Weil bound states that if $C/\F_q$ is a smooth projective geometrically integral curve of genus one, then
\[
  \bigl|\#C(\F_q)-(q+1)\bigr|\le 2\sqrt q,
\]
see Stichtenoth~\cite[Chapter~5, \S5.2]{Stichtenoth2009}.  In particular, for an elliptic curve $E/\F_q$ with Frobenius trace $t_E=q+1-\#E(\F_q)$, one has $|t_E|\le2\sqrt q$.  For an elliptic curve $E/\Q$, $j(E)$ denotes its $j$-invariant.  We say that $E$ has complex multiplication (CM) if $\operatorname{End}_{\overline{\Q}}(E)$ is strictly larger than $\Z$.  The CM facts used later are cited where they enter.

The local two-target formulas used here are proved in~\cite{FactCrypt}.  We study the following related problem, which is relevant to the coefficient restrictions in deterministic factorization.  We seek to choose a primitive root $\gamma$ modulo $\ell$ and use the two coefficients $\gamma$ and $\gamma^2$.  Since $\gamma$ is a nonsquare and $\gamma^2$ is a square, this pair is the smallest model that simultaneously represents both square classes of coefficients.  The desired local inequality is $|H_\gamma(N;\ell)\cup H_{\gamma^2}(N;\ell)|\le 3\ell/4+1$ for every $N\in\Fps$.  The exact formula in~\cite{FactCrypt} shows that this pointwise bound is equivalent to a bounded Frobenius-trace condition on the elliptic curve $E_\gamma:Y^2=(X^2-1)(X^2-\gamma)$.

We consider two local constructions and their use in deterministic factorization.  The first is trace selection, which asks for primes $\ell$ for which there exists a primitive root $\gamma$ whose Frobenius trace lies in the required short interval.  The second uses the residue-known Harvey--Hittmeir lattice directly.  Rather than imposing an arbitrary CRT class on a Lehman approximant, we keep the coefficient lattice of~\cite[Proposition~3.3]{HarveyHittmeir2022}, compute one further local digit of the actual babystep index, and show that this digit is a quadratic image modulo every usable prime divisor of the auxiliary modulus.  A candidate-dependent translation moves a fixed quarter of the local residue space into the complement of that image.  Combining these local restrictions with Gauss reduction, sparse multipoint evaluation, and the unrestricted large-order theorem of Harvey and Hittmeir~\cite{HarveyHittmeir2026} gives the factorization bound stated in Theorem~\ref{thm:factorization-main}.

\section{The trace-selection approach}
\label{sec:trace-selection}

We recall the local statements from~\cite{FactCrypt} in the notation used here.  First, $|H_a(N;\ell)|=(\ell+\chi(aN))/2$, as proved in~\cite[Lemma~3.2]{FactCrypt}.  For two distinct coefficients $a,b\in\Fps$, put $\U_{a,b}(N;\ell)=H_a(N;\ell)\cup H_b(N;\ell)$, $A=4aN$, and $B=4bN$, and set $r_A=1+\chi(A)$, $r_B=1+\chi(B)$.  The exact formula proved in~\cite[Theorem~3.9]{FactCrypt} is
\[
  4|\Fp\setminus\U_{a,b}(N;\ell)|
  =\ell+2-r_A-r_B+K(A,B)+r_B\chi(B-A)+r_A\chi(A-B),
\]
where $K(A,B)=\sum_{z\in\Fp}\chi((z^2-A)(z^2-B))$.  Since $Y^2=(X^2-A)(X^2-B)$ has a smooth projective genus-one model for $A,B\ne0$ and $A\ne B$, the genus-one Hasse--Weil bound stated above yields the following uniform consequence through the exact formula of~\cite[Theorem~3.9]{FactCrypt}.

\begin{corollary}
\label{cor:hasse-union-new}
For any prime $\ell>2$, any $N\in\Fps$, and any distinct $a,b\in\Fps$,
\[
  \left||H_a(N;\ell)\cup H_b(N;\ell)|-\frac{3\ell}{4}\right|
  \le \frac{2\sqrt\ell+7}{4}.
\]
\end{corollary}

Indeed, the quartic has two rational points at infinity, so $K(A,B)=-1-t$ with $|t|\le2\sqrt\ell$.  The remaining terms in the exact complement formula have total absolute value at most $6$, which gives the stated bound.  The primitive-root specialization is sharper because the two possible twists can be written in terms of a single trace.

\begin{theorem}
\label{thm:primitive-root-two-values-new}
Let $\gamma$ be a primitive root modulo $\ell$, set $A=4\gamma N$, and let $t_\gamma=\ell+1-\#E_\gamma(\Fp)$ be the Frobenius trace of
\[
  E_\gamma: Y^2=(X^2-1)(X^2-\gamma).
\]
Here $E_\gamma$ denotes the smooth projective model, including its two rational points at infinity.  Then
\[
  |H_\gamma(N;\ell)\cup H_{\gamma^2}(N;\ell)|
  =
  \begin{cases}
  \displaystyle\frac{3\ell+1+t_\gamma-2\chi(1-\gamma)}{4},& A\text{ is a square},\\[.8em]
  \displaystyle\frac{3\ell+1-t_\gamma+2\chi(\gamma-1)}{4},& A\text{ is a nonsquare}.
  \end{cases}
\]
\end{theorem}

\begin{proof}
This is the primitive-root specialization of~\cite[Theorem~3.13]{FactCrypt}.  We include the character-sum calculation to fix the sign convention.  Put $K_\gamma=K(1,\gamma)$.  Counting the two points at infinity gives $K_\gamma=-1-t_\gamma$.  For $A\ne0$, substitute $u=z^2/A$ in the character sum.  Each $u$ has $1+\chi(A)\chi(u)$ preimages, and therefore
\[
 K(A,\gamma A)
 =\sum_u\bigl(1+\chi(A)\chi(u)\bigr)
       \chi((u-1)(u-\gamma))
 =-1+\chi(A)(K_\gamma+1)
 =-1-\chi(A)t_\gamma.
\]
Because $\gamma$ is a nonsquare, exactly one of $A$ and $\gamma A$ is a square.  Thus $r_A+r_{\gamma A}=2$.  If $A$ is a square, the boundary contribution in the complement formula is $2\chi(1-\gamma)$.  If $A$ is a nonsquare, it is $-2\chi(\gamma-1)$.  Substitution gives the two formulas.
\end{proof}

Thus, for fixed $\ell$ and $\gamma$, the primitive-root two-target union has at most two possible sizes as $N$ varies.  The exact pointwise criterion is the following.

\begin{corollary}
\label{cor:primitive-root-trace-criterion}
Let $\gamma$ be a primitive root modulo $\ell$, and let $t_\gamma$ be the Frobenius trace of the smooth projective curve $E_\gamma$ given by $Y^2=(X^2-1)(X^2-\gamma)$.  Then
\[
  |H_\gamma(N;\ell)\cup H_{\gamma^2}(N;\ell)|\le \frac{3\ell}{4}+1
  \quad\text{for every }N\in\Fps
\]
if and only if
\[
  -3+2\chi(\gamma-1)\le t_\gamma\le 3+2\chi(1-\gamma).
\]
\end{corollary}

\begin{proof}
The two square classes of $A=4\gamma N$ both occur as $N$ varies over $\Fps$.  The first formula in Theorem~\ref{thm:primitive-root-two-values-new} is at most $3\ell/4+1$ exactly when $t_\gamma\le3+2\chi(1-\gamma)$.  The second is at most that bound exactly when $t_\gamma\ge-3+2\chi(\gamma-1)$.  Both inequalities are therefore necessary and sufficient.
\end{proof}

The interval in Corollary~\ref{cor:primitive-root-trace-criterion} is much shorter than the Hasse--Weil interval $|t_\gamma|\le2\sqrt\ell$.  Therefore the problem is not to improve Hasse--Weil uniformly.  It is to find primitive roots whose traces lie in the specified bounded interval.

\begin{definition}
\label{def:admissible-primitive-root}
A primitive root $\gamma\in\Fps$ is called \emph{trace-admissible for $\ell$} if, for $E_\gamma:Y^2=(X^2-1)(X^2-\gamma)$ and $t_\gamma=\ell+1-\#E_\gamma(\Fp)$, one has
\[
  -3+2\chi(\gamma-1)\le t_\gamma\le 3+2\chi(1-\gamma).
\]
The prime $\ell$ is called \emph{trace-admissible} if it has at least one trace-admissible primitive root.
\end{definition}

The trace-admissibility interval has bounded length, but its length depends on the congruence class of $\ell$ modulo $4$ and on the sign of $\chi(\gamma-1)$.  This elementary point is useful when estimating how rare trace-admissible primitive roots should be.

\begin{lemma}
\label{lem:admissibility-interval-length}
Let $\ell>2$ be prime and let $\gamma\in\Fps$ with $\gamma\ne1$.  Put
\[
  I_\ell(\gamma)=\{t\in\Z:-3+2\chi(\gamma-1)\le t\le 3+2\chi(1-\gamma)\}.
\]
Then
\[
  |I_\ell(\gamma)|=7+2\chi(1-\gamma)-2\chi(\gamma-1).
\]
Consequently, if $\ell\equiv1\pmod4$, then $|I_\ell(\gamma)|=7$, while if $\ell\equiv3\pmod4$, then $|I_\ell(\gamma)|=7-4\chi(\gamma-1)$, hence it is either $3$ or $11$.  In particular, the interval contains $3$, $7$, or $11$ integers.
\end{lemma}

\begin{proof}
The number of integers in the interval $[A,B]$ is $B-A+1$, with $A=-3+2\chi(\gamma-1)$ and $B=3+2\chi(1-\gamma)$.  This gives the displayed formula.  Since $1-\gamma=-(\gamma-1)$, one has $\chi(1-\gamma)=\chi(-1)\chi(\gamma-1)$.  The two cases follow from $\chi(-1)=1$ for $\ell\equiv1\pmod4$ and $\chi(-1)=-1$ for $\ell\equiv3\pmod4$.
\end{proof}

\begin{proposition}
\label{prop:admissibility-equivalence}
Let $\ell>2$ be prime.  A primitive root $\gamma\in\Fps$ is trace-admissible if and only if the two-target local inequality
\[
  |H_\gamma(N;\ell)\cup H_{\gamma^2}(N;\ell)|\le \frac{3\ell}{4}+1
\]
holds for every $N\in\Fps$.
\end{proposition}

\begin{proof}
This is Corollary~\ref{cor:primitive-root-trace-criterion} written in the terminology of Definition~\ref{def:admissible-primitive-root}.
\end{proof}

\begin{proposition}
\label{prop:certificate}
Let $\ell>2$ be prime and let $\gamma\in\Fps$.  To certify that $\gamma$ is a trace-admissible primitive root, provide a complete prime factorization of $\ell-1$, verify $\gamma^{(\ell-1)/r}\ne1$ for every prime divisor $r$ of $\ell-1$, and compute the exact point count $\#E_\gamma(\Fp)$.  The resulting trace must satisfy the interval in Definition~\textup{\ref{def:admissible-primitive-root}}.
\end{proposition}

\begin{proof}
The first two checks prove that $\gamma$ generates the cyclic group $\Fps$.  The point count determines $t_\gamma$ exactly.  The remaining interval check is the definition of trace-admissibility.
\end{proof}

The preceding proposition gives a local certificate.  A factorization speedup also requires compatibility with the global coefficient search.  For a CRT modulus $m=\prod_i\ell_i$, one can nevertheless record a product certificate for the local sets.  If each $\ell_i$ has a trace-admissible primitive root $\gamma_i$, then for every residue $N\in(\Z/m\Z)^*$ one has
\[
  \prod_i |H_{\gamma_i}(N;\ell_i)\cup H_{\gamma_i^2}(N;\ell_i)|
  \le
  \prod_i \left(\frac{3\ell_i}{4}+1\right).
\]
For $\ell_i>4$, each displayed factor is strictly less than $\ell_i$, so the product is smaller than the full CRT residue space when all $\ell_i\ge5$.  To use this reduction in the factorization algorithm, one must also show that the corresponding coefficient ratios are compatible with the global Lehman search.

\begin{problem}
\label{prob:trace-selection}
Determine whether there are infinitely many trace-admissible primes.  More quantitatively, estimate the number of primes $\ell\le X$ for which there exists a primitive root $\gamma$ modulo $\ell$ such that the curve $E_\gamma:Y^2=(X^2-1)(X^2-\gamma)$ has trace $t_\gamma$ in the interval
\[
  -3+2\chi(\gamma-1)\le t_\gamma\le 3+2\chi(1-\gamma).
\]
\end{problem}

\begin{remark}
\label{rem:distribution-warning}
The trace-selection problem is not settled by a direct appeal to Sato--Tate heuristics.  The required statement concerns the traces of the family $E_\gamma$ when $\gamma$ ranges only over primitive roots modulo a varying prime $\ell$.  A proof must therefore combine quantitative trace distribution for the family with nontrivial control of the primitive-root restriction.  Neither the Hasse--Weil bound nor the exact union formula quoted above supplies such a result.  A fixed-generator argument also requires care.  Fixing an integer $g$ and asking simultaneously that $g\bmod\ell$ be a primitive root and that the fixed curve $E_g$ have trace in a bounded interval is a much thinner problem than asking whether some primitive root $\gamma$ modulo $\ell$ has the required trace.
\end{remark}

\begin{remark}
\label{rem:expected-number}
For primes $\ell\not\equiv5\pmod8$, a heuristic count suggests the scale $\varphi(\ell-1)/\sqrt\ell$.  There are $\varphi(\ell-1)$ primitive roots, while the trace ranges over $O(\sqrt\ell)$ integers and the admissibility interval contains $O(1)$ integers.  This estimate presumes that admissible trace values occur with the expected frequency under the primitive-root restriction.  Congruences can exclude all such values, as they do for $\ell\equiv5\pmod8$, so a continuous Sato--Tate density does not justify a leading constant at this scale.  The estimate is only a heuristic for the three remaining residue classes.
\end{remark}

The reported primitive-root computations suggest a classification by residue class modulo $8$.  Checking all primes up to $10^4$ suggests that all primes congruent to $1$ or $3$ modulo $8$ are trace-admissible, that no prime congruent to $5$ modulo $8$ is trace-admissible, and that in the class $7\bmod 8$ the only exceptions are $\ell=7$ and $\ell=823$.  After omitting the class $5\bmod8$, the same pattern in the classes $1,3,7\bmod8$ has been checked up to $10^5$.  The negative assertion in the class $5\bmod8$ is not merely experimental because it follows from the mod-$8$ congruence proved below.  The positive assertions remain conjectural, except for the GRH-conditional infinite subfamily in Theorem~\ref{thm:grh-cm-trace-admissible} below.

\begin{conjecture}
\label{conj:mod-eight-pattern}
Let $\ell>2$ be prime.  The proposed classification is as follows.
\begin{enumerate}[label=\textup{(\alph*)}]
\item If $\ell\equiv1$ or $3\pmod8$, then $\ell$ is trace-admissible.
\item If $\ell\equiv5\pmod8$, then $\ell$ is not trace-admissible, as proved in Proposition~\ref{prop:five-mod-eight-congruence}.
\item If $\ell\equiv7\pmod8$, then $\ell$ is trace-admissible except for $\ell=7$ and $\ell=823$.
\end{enumerate}
\end{conjecture}

The next proposition proves the second part of Conjecture~\ref{conj:mod-eight-pattern}.  The proof uses only the elementary point count on the quartic model.  When $\ell\equiv5\pmod8$ and $\gamma$ is a primitive root, the affine points with $y\ne0$ occur in orbits of size four under $(x,y)\mapsto(-x,y)$ and $(x,y)\mapsto(x,-y)$, and the parity of the number of such orbits is determined by a product over the nonzero squares of $\F_\ell$.

\begin{proposition}
\label{prop:five-mod-eight-congruence}
Let $\ell\equiv1\pmod4$ be prime and let $\gamma$ be any nonsquare modulo $\ell$.  Let
\[
  E_\gamma:Y^2=(X^2-1)(X^2-\gamma),
  \qquad
  t_\gamma=\ell+1-\#E_\gamma(\F_\ell),
\]
where the smooth projective model of the quartic has the two rational points at infinity.  Then
\[
  t_\gamma\equiv \ell-1+2\chi(\gamma-1)\pmod8.
\]
Consequently no prime $\ell\equiv5\pmod8$ is trace-admissible.  If $\ell\equiv1\pmod8$ and $\gamma$ is primitive, then $\gamma$ is trace-admissible if and only if $t_\gamma=2\chi(\gamma-1)$.
\end{proposition}

\begin{proof}
Since $\ell\equiv1\pmod4$, one has $\chi(-1)=1$.  Write $h=\chi(2)$ and $e=\chi(\gamma-1)$.  The parameter $\gamma$ is a nonsquare by assumption.  The polynomial $f(X)=(X^2-1)(X^2-\gamma)$ has exactly two roots in $\F_\ell$, namely $X=\pm1$, because $\gamma$ is a nonsquare.  Also $f(0)=\gamma$ is a nonsquare, so there are no affine points above $X=0$.  Apart from the two points $(\pm1,0)$ and the two points at infinity, all rational points have $x\ne0$, $x\ne\pm1$, and $y\ne0$, and therefore occur in orbits of size four under $(x,y)\mapsto(-x,y)$ and $(x,y)\mapsto(x,-y)$.  Let $M$ be the number of orbits of this kind.  Then $\#E_\gamma(\F_\ell)=4+4M$.  Thus it remains to determine $M$ modulo $2$.

Let $Q\subset\F_\ell^*$ be the set of nonzero squares and let $S=Q\setminus\{1\}$.  The orbits counted by $M$ are in bijection with those $u=x^2\in S$ for which $(u-1)(u-\gamma)$ is a nonzero square.  Hence $M=\#\{u\in S:\chi((u-1)(u-\gamma))=1\}$.  Since $|Q|=(\ell-1)/2\equiv0\pmod2$, the set $S$ has odd cardinality.  Therefore, if $P=\prod_{u\in S}\chi((u-1)(u-\gamma))$, then $(-1)^M=-P$.  We compute $P$.  Write $n=(\ell-1)/2$.  The polynomial with roots $Q$ is $X^n-1$.  Differentiating at $X=1$ gives $\prod_{u\in Q,\,u\ne1}(1-u)=n$, hence $\prod_{u\in S}(u-1)=(-1)^{n-1}n=-n$.  Because $n\equiv-1/2\pmod\ell$, and because $\chi(-1)=1$, we obtain $\chi(\prod_{u\in S}(u-1))=h$.  Next, since $\gamma$ is a nonsquare, $\gamma^n=-1$.  As $n$ is even, $\prod_{u\in Q}(u-\gamma)=\prod_{u\in Q}(\gamma-u)=\gamma^n-1=-2$.  Removing the factor $u=1$ gives $\prod_{u\in S}(u-\gamma)=-2/(1-\gamma)$, and therefore $\chi(\prod_{u\in S}(u-\gamma))=\chi(-2)\chi(1-\gamma)=he$, because $\chi(-2)=h$ and $\chi(1-\gamma)=\chi(\gamma-1)$ in the class $\ell\equiv1\pmod4$.  Hence $P=e$ and $(-1)^M=-e$.  It follows that $M$ is odd if $e=1$ and even if $e=-1$.  Consequently
\[
  \#E_\gamma(\F_\ell)\equiv
  \begin{cases}
  0\pmod8,& \chi(\gamma-1)=1,\\
  4\pmod8,& \chi(\gamma-1)=-1.
  \end{cases}
\]
Subtracting these point counts from $\ell+1$ gives
\[
  t_\gamma=\ell+1-\#E_\gamma(\F_\ell)
  \equiv \ell-1+2\chi(\gamma-1)\pmod8.
\]

Finally, for $\ell\equiv5\pmod8$ we have $\chi(1-\gamma)=\chi(\gamma-1)$.  If $\chi(\gamma-1)=1$, the trace-admissibility interval is $[-1,5]$, but the congruence gives $t_\gamma\equiv6\pmod8$, and no integer in $[-1,5]$ is congruent to $6$ modulo $8$.  If $\chi(\gamma-1)=-1$, the interval is $[-5,1]$, but the congruence gives $t_\gamma\equiv2\pmod8$, and no integer in $[-5,1]$ is congruent to $2$ modulo $8$.  Thus no primitive root is trace-admissible, and $\ell$ is not trace-admissible.  For $\ell\equiv1\pmod8$, the same congruence is $t_\gamma\equiv2e\pmod8$.  The interval is $[-1,5]$ if $e=1$ and $[-5,1]$ if $e=-1$, and its unique integer in the required congruence class is $2e$.  This proves the final assertion.
\end{proof}

\begin{problem}
\label{prob:positive-mod-eight-classes}
Prove the remaining positive assertions suggested by computation.  Every prime $\ell\equiv1$ or $3\pmod8$ should be trace-admissible, and every prime $\ell\equiv7\pmod8$ should be trace-admissible except $7$ and $823$.  The class $3\pmod8$ contains the GRH-conditional infinite subfamily supplied by Theorem~\ref{thm:grh-cm-trace-admissible}, and the computational evidence suggests that trace admissibility holds throughout these classes, with the stated exceptions, without fixing $\gamma=2$.
\end{problem}

\begin{problem}
\label{prob:seven-mod-eight-exceptions}
Explain why $\ell=823$ is non-admissible in the class $7\bmod8$, and prove that no further exceptions occur.  For $\ell=7$, the primitive roots are $\gamma=3$ and $\gamma=5$.  Since $7\equiv3\pmod4$, one has $\chi(1-\gamma)=-\chi(\gamma-1)$.  For both primitive roots, $\chi(\gamma-1)=1$, so the admissibility interval is $[-1,1]$.  Direct point counting gives $t_3=4$ and $t_5=-4$, and neither trace lies in that interval.  For $\ell=823$, direct computation confirms that all primitive roots have traces outside their admissibility intervals.  The problem is to identify the arithmetic reason for this isolated failure.
\end{problem}

\begin{hypothesis}
\label{hyp:hooley-artin-ap}
For every squarefree integer $d\ge1$, let
\[
 K_d=\Q(\zeta_{\mathrm{lcm}(4,d)},2^{1/d}).
\]
The Dedekind zeta function of each $K_d$ satisfies the Generalized Riemann Hypothesis.  This is a sufficient GRH assumption for the arithmetic-progression form of Hooley's argument used below.
\end{hypothesis}

\begin{lemma} 
\label{lem:hooley-2-mod4}
Assume Hypothesis~\textup{\ref{hyp:hooley-artin-ap}}.  Then there is a constant $A_{2,3(4)}>0$ such that
\[
  \#\{\ell\le X:\ell\equiv3\pmod4\text{ and }2\text{ is a primitive root modulo }\ell\}
  \sim A_{2,3(4)}\operatorname{Li}(X).
\]
In particular, the set of such primes has positive density among all primes.
\end{lemma}

\begin{proof}
This is the arithmetic-progression version of the GRH-conditional Artin primitive-root theorem.  Hooley's original theorem treats primitive roots without a prescribed congruence class~\cite{Hooley1967}, and the extension to primes in arithmetic progressions is contained in Lenstra's GRH-conditional formulation~\cite{Lenstra1977}, with explicit density formulas discussed by Moree~\cite{Moree2008}.  Specializing those results to the integer $2$ and the progression $3\pmod4$ gives the asserted asymptotic, under Hypothesis~\ref{hyp:hooley-artin-ap}.

For completeness, the relevant inclusion--exclusion starts from the primitive-root indicator $\sum_{d\mid \ell-1}\mu(d)\,\mathbf 1[2^{(\ell-1)/d}\equiv1\pmod\ell]$.  For each squarefree $d$, the simultaneous conditions $2^{(\ell-1)/d}\equiv1\pmod\ell$ and $\ell\equiv3\pmod4$ are encoded by Frobenius classes in $\Q(\zeta_{\mathrm{lcm}(4,d)},2^{1/d})$.  The GRH hypothesis stated above supplies the analytic input needed to apply effective Chebotarev uniformly over these fields and to control the inclusion--exclusion tail.  Moree's density formula gives, in this specialization,
\[
 A_{2,3(4)}=\frac12\prod_{r\ {\rm prime}}
       \left(1-\frac1{r(r-1)}\right)>0.
\]
The product is positive because each factor is positive and the sum of $1/(r(r-1))$ over the primes converges.  Thus positivity follows from the density formula itself.
\end{proof}

Before proving the GRH-conditional infinitude theorem, we make the special CM input completely explicit.  The computation of the $j$-invariant uses the classical invariants of a binary quartic.  For $f(X)=aX^4+bX^3+cX^2+dX+e$ with nonzero discriminant, put $I=12ae-3bd+c^2$ and $J=72ace+9bcd-27ad^2-27b^2e-2c^3$.  The Jacobian elliptic curve of $Y^2=f(X)$ has $j=6912I^3/(4I^3-J^2)$.  We use this formula only in the explicit quartic below, where the substitution is elementary.

\begin{lemma}
\label{lem:j-invariant-family}
For $\gamma\ne0,1$ over a field of characteristic different from $2$, the smooth projective model of
\[
  E_\gamma:\quad Y^2=(X^2-1)(X^2-\gamma)
\]
has
\[
  j(E_\gamma)=16\frac{(\gamma^2+14\gamma+1)^3}{\gamma(\gamma-1)^4}.
\]
In particular, $j(E_2)=287496=66^3$.
\end{lemma}

\begin{proof}
Here $f(X)=X^4-(1+\gamma)X^2+\gamma$, so $a=1$, $b=d=0$, $c=-(1+\gamma)$, and $e=\gamma$.  The binary-quartic invariants are therefore $I=12\gamma+(1+\gamma)^2=\gamma^2+14\gamma+1$ and $J=-72\gamma(1+\gamma)+2(1+\gamma)^3=2(1+\gamma)(\gamma^2-34\gamma+1)$.  A direct expansion gives $4I^3-J^2=432\gamma(\gamma-1)^4$.  Substitution into $j=6912I^3/(4I^3-J^2)$ gives
\[
  j(E_\gamma)=16\frac{(\gamma^2+14\gamma+1)^3}{\gamma(\gamma-1)^4}.
\]
The identity is first obtained over $\Q(\gamma)$ and extends to every odd characteristic, including characteristic $3$, because its final denominator is nonzero for $\gamma\ne0,1$.  Putting $\gamma=2$ gives $I=33$ and therefore $j(E_2)=16\cdot33^3/2=8\cdot33^3=66^3=287496$.
\end{proof}

\begin{lemma}
\label{lem:cm-trace-zero}
Let $\ell>2$ be a prime with $\ell\equiv3\pmod4$.  For the smooth projective curve $E_2$ given by $Y^2=(X^2-1)(X^2-2)$ over $\F_\ell$, the Frobenius trace is $t_2(\ell)=0$.
\end{lemma}

\begin{proof}
Over $\Q$, Lemma~\ref{lem:j-invariant-family} gives $j(E_2)=287496$.  This is the singular modulus attached to the order of discriminant $-16$ in $\Q(i)$, as recorded in Cox~\cite[Chapter~11]{Cox1989}, where the class polynomials and singular moduli for imaginary quadratic orders are described, or Silverman~\cite[Chapter~II, Section~6]{Silverman1994}.  Hence $E_2$ has complex multiplication by an order in $K=\Q(i)$.

Since $\ell\equiv3\pmod4$, the prime $\ell$ is inert in $K$.  Deuring's reduction theorem for elliptic curves with complex multiplication states that, at a prime of good reduction inert in the CM field, the reduction is supersingular, as stated in Silverman~\cite[Chapter~II, Section~6]{Silverman1994}.  The quartic defining $E_2$ is squarefree modulo every odd prime and has a rational point at infinity, so $E_2$ has good reduction at every $\ell>2$.  Therefore the reduction of $E_2$ modulo $\ell$ is supersingular.  For a supersingular elliptic curve over $\F_\ell$ with $\ell>3$, the trace is divisible by $\ell$.  Combined with Hasse's bound $|t_2(\ell)|\le2\sqrt\ell<\ell$, this gives $t_2(\ell)=0$.  For the remaining prime $\ell=3$, direct counting gives
\[
  \#E_2(\F_3)=4=3+1,
\]
so $t_2(3)=0$ as well.
\end{proof}

Combining the trace-zero specialization with the primitive-root density gives the following conditional infinitude result.
\begin{theorem}
\label{thm:grh-cm-trace-admissible}
Assume Hypothesis~\textup{\ref{hyp:hooley-artin-ap}}.  Then
\[
  \#\{\ell\le X:\ell\text{ is trace-admissible}\}
  \gg \frac{X}{\log X}.
\]
More precisely, every prime $\ell>2$ such that $\ell\equiv3\pmod4$ and $2$ is a primitive root modulo $\ell$ is trace-admissible.
\end{theorem}

\begin{proof}
Let $\ell>2$ satisfy $\ell\equiv3\pmod4$, and suppose that $2$ is a primitive root modulo $\ell$.  Put $\gamma=2\in\F_\ell^{\!*}$.  Then $\gamma$ is a primitive root and $t_\gamma=0$ by Lemma~\ref{lem:cm-trace-zero}.  Moreover
\[
  \chi(\gamma-1)=\chi(1)=1,
  \qquad
  \chi(1-\gamma)=\chi(-1)=-1,
\]
because $\ell\equiv3\pmod4$.  Therefore the trace-admissibility interval is
\[
  -3+2\chi(\gamma-1)\le t\le 3+2\chi(1-\gamma),
  \qquad\text{i.e.}\qquad
  -1\le t\le1.
\]
Since $t_\gamma=0$, the primitive root $\gamma=2$ is trace-admissible.  The lower bound now follows from Lemma~\ref{lem:hooley-2-mod4}.
\end{proof}

\begin{remark}
\label{rem:scope-of-grh-result}
Theorem~\ref{thm:grh-cm-trace-admissible} proves infinitude only for the specialization $\gamma=2$.  This gives a positive-density conditional subfamily, but it does not account for the broader computational pattern in Conjecture~\ref{conj:mod-eight-pattern}.  In particular, it gives no unconditional information about the classes $\ell\equiv1\pmod8$ and $\ell\equiv7\pmod8$, and it should not be interpreted as a replacement for the finite-field trace-selection problem.
\end{remark}

Corollary~\ref{cor:hasse-union-new} gives an unconditional alternative if the global argument needs only an asymptotic local saving instead of the pointwise bound $3\ell/4+1$.  It yields a product estimate for arbitrary distinct coefficient pairs.

\begin{proposition}
\label{prop:unconditional-crt-product}
Let $\ell_1,\ldots,\ell_s$ be distinct odd primes, let $m=\prod_i\ell_i$, and for each $i$ choose distinct coefficients $a_i,b_i\in\F_{\ell_i}^{\!*}$.  For every $N\in(\Z/m\Z)^*$ one has
\[
  \prod_{i=1}^s |H_{a_i}(N;\ell_i)\cup H_{b_i}(N;\ell_i)|
  \le
  \prod_{i=1}^s\left(\frac{3\ell_i}{4}+\frac{2\sqrt{\ell_i}+7}{4}\right).
\]
Equivalently, the local density at $\ell_i$ is at most
\[
  \frac34+O(\ell_i^{-1/2}).
\]
\end{proposition}

\begin{proof}
Apply Corollary~\ref{cor:hasse-union-new} separately at each prime $\ell_i$ and multiply the resulting inequalities.
\end{proof}

This estimate is weaker than trace admissibility at a fixed small prime, but it is unconditional and uniform in the coefficient pair.  One can therefore choose the auxiliary primes large enough that the $O(\ell_i^{-1/2})$ term is negligible and then study the coefficient-selection problem for arbitrary two-target CRT classes.  This avoids primitive-root trace selection, at the cost of using an asymptotic local bound instead of an exact one.

The next hypothesis states a sufficient finite-field estimate for trace selection in the three residue classes not excluded above.  It is formulated at a fixed prime $\ell$ and does not follow formally from GRH for number-field $L$-functions.

\begin{hypothesis}
\label{hyp:mixed-trace-primitive-root}
There are constants $c_0>0$, $C_0>0$, $\eta>0$, and $\ell_0$ such that, for every prime $\ell\ge\ell_0$ with $\ell\not\equiv5\pmod8$,
\[
  \#\{\gamma\in\F_{\ell}^{\!*}:\gamma\text{ is a primitive root and }t_\gamma\in I_\ell(\gamma)\}
  \ge c_0\frac{\varphi(\ell-1)}{\sqrt{\ell}}-C_0\ell^{1/2-\eta},
\]
where $I_\ell(\gamma)$ is the interval of Lemma~\textup{\ref{lem:admissibility-interval-length}}.
\end{hypothesis}

\begin{proposition}
\label{prop:mixed-estimate-implies-admissible}
Under Hypothesis~\textup{\ref{hyp:mixed-trace-primitive-root}}, all sufficiently large primes $\ell\not\equiv5\pmod8$ are trace-admissible.
\end{proposition}

\begin{proof}
The main term is $\gg \varphi(\ell-1)/\sqrt{\ell}$.  Since $\varphi(n)\gg n/\log\log n$ for $n\ge3$, this is $\gg \sqrt{\ell}/\log\log\ell$, which dominates the error term $O(\ell^{1/2-\eta})$ for sufficiently large $\ell$.  Hence the displayed count is positive for all sufficiently large $\ell$ in these three residue classes.
\end{proof}

Hypothesis~\ref{hyp:mixed-trace-primitive-root} does not follow from the GRH assumptions used in this paper.  For a fixed integer $g$, the condition that $g$ be a primitive root is Artin-type, whereas for non-CM $E_g$ the additional demand of a trace in a bounded interval is of Lang--Trotter type.  CM curves can instead have trace zero on a positive-density set of primes, as in the specialization $g=2$ above.  Hypothesis~\ref{hyp:mixed-trace-primitive-root} is a different finite-field statement in which $\gamma$ varies with $\ell$, and one needs cancellation in mixed sums over primitive roots inside the family $E_\gamma$.

The following observation is useful for sharpening the target, but it also explains why the primitive-root condition cannot be disposed of by a bare Plancherel argument.  Let
\[
  A_\ell=\{\gamma\in\F_\ell^*\setminus\{1\}: t_\gamma\in I_\ell(\gamma)\}.
\]
If one writes the primitive-root indicator on the cyclic group $\F_\ell^*$ as
\[
  \mathbf 1_{\mathrm{prim}}(\gamma)
  =\sum_{\psi\in\widehat{\F_\ell^*}} c(\psi)\psi(\gamma),
\]
then the principal coefficient is $c(\chi_0)=\varphi(\ell-1)/(\ell-1)$, but the remaining coefficients are not small enough individually to force positivity from the size of $A_\ell$ alone.  More explicitly, the standard formula on a cyclic group of order $n=\ell-1$ is
\[
  \mathbf 1_{\mathrm{prim}}(x)
  =\frac{\varphi(n)}{n}
   \sum_{d\mid n}\frac{\mu(d)}{\varphi(d)}
   \sum_{\substack{\psi\in\widehat{\F_\ell^*}\\ \operatorname{ord}(\psi)=d}}\psi(x).
\]
Thus a character of order $d$ has coefficient $(\varphi(n)/n)\mu(d)/\varphi(d)$.  Applying the Plancherel identity stated above to $f=\mathbf 1_{A_\ell}$ gives
\[
  \sum_{\psi}|\sum_{\gamma\in A_\ell}\psi(\gamma)|^2=(\ell-1)|A_\ell|,
\]
and applying the Cauchy--Schwarz inequality stated above to the weighted character sum, the nonprincipal contribution is bounded by a quantity of order
\[
  \left(\frac{\varphi(n)}{n}\right)
  \left(\sum_{d\mid n}\frac{\mu(d)^2}{\varphi(d)}\right)^{1/2}
  ((\ell-1)|A_\ell|)^{1/2},
\]
which is too large when $|A_\ell|$ is only of order $\sqrt\ell$.  Therefore this approach to a density-one theorem within the classes $1,3,7\pmod8$ would require estimates stronger than the resulting Cauchy--Schwarz bound for the twisted trace sums
\[
  \sum_{\gamma\in A_\ell}\psi(\gamma),\qquad \psi\ne\chi_0,
\]
not merely Plancherel.  Hypothesis~\ref{hyp:mixed-trace-primitive-root} records a sufficient lower bound for the primitive-root count.  It is restricted to $\ell\not\equiv5\pmod8$ because Proposition~\ref{prop:five-mod-eight-congruence} rules out the omitted class.  In particular, trace-admissible primes have upper relative density at most $3/4$ among all primes.

\section{The Harvey--Hittmeir lattice and quadratic-image sieving}
\label{sec:hh-gap-interface}

We now consider the global search.  We keep the residue-known coefficient
lattice of Harvey and Hittmeir~\cite[Proposition~3.3]{HarveyHittmeir2022}
and modify only the choice of a short lattice vector and the final collision
search.  Modulo each usable prime divisor of the auxiliary modulus, the
Harvey--Hittmeir babystep index is the image of a quadratic polynomial.  A
giantstep-dependent translation places a fixed set containing about one
quarter of the residues outside that image.  Restricting to sufficiently
large prime divisors of the auxiliary modulus allows us to choose a nearby
reduced-lattice vector for which all these local quadratic polynomials are
nondegenerate.

Let $N=pq$, let $m\ge1$ with $\gcd(m,N)=1$, and let
$\sigma\in(\Z/m\Z)^*$ satisfy $p\equiv\sigma\pmod m$.  Fractions whose
denominators are coprime to $m$ are interpreted in the localization
$\Z_{(m)}$, and reduction modulo a divisor of a power of $m$ is performed by
using the corresponding modular inverse.

\begin{lemma}
\label{lem:hh-residue-lift}
Suppose that $a,b\in\Z$ satisfy
\begin{equation}
\label{eq:hh-lattice-congruence}
  -a\frac{N}{\sigma}+b\sigma\equiv0\pmod m.
\end{equation}
Then
\begin{equation}
\label{eq:hh-modm2}
  aq+bp\equiv a\frac{N}{\sigma}+b\sigma\pmod{m^2}.
\end{equation}
Equivalently, the solutions of~\eqref{eq:hh-lattice-congruence} form the
rank-two lattice
\begin{equation*}
%\label{eq:hh-lattice}
  L_\sigma
  =\{(a,b)\in\Z^2:b\equiv\eta a\pmod m\},
  \qquad
  \eta\equiv N\sigma^{-2}\pmod m,
\end{equation*}
which has basis $(1,\eta),(0,m)$ and determinant $m$.
\end{lemma}

\begin{proof}
Write $p=\sigma(1+t)$ with $t=(p-\sigma)/\sigma$.  Then $t\in\Z_{(m)}$ and
$t\equiv0\pmod m$.  Hence
\[
  q=\frac N\sigma(1+t)^{-1}
  \equiv\frac N\sigma(1-t)\pmod{m^2}.
\]
Therefore
\[
  aq+bp
  \equiv
  a\frac N\sigma+b\sigma
  +t\left(-a\frac N\sigma+b\sigma\right)
  \pmod{m^2}.
\]
The last term is divisible by $m^2$ under
\eqref{eq:hh-lattice-congruence}, proving~\eqref{eq:hh-modm2}.  Multiplying
\eqref{eq:hh-lattice-congruence} by the unit $\sigma$ modulo $m$ yields
$b\sigma^2\equiv aN\pmod m$, which is equivalent to the displayed lattice
description.  The basis and determinant are immediate.
\end{proof}

\begin{lemma}
\label{lem:local-divisibility-dichotomy}
Assume that $m$ is squarefree and $(a,b)\in L_\sigma$.  For every prime
$\ell\mid m$,
\[
  \ell\mid a\quad\Longleftrightarrow\quad\ell\mid b.
\]
\end{lemma}

\begin{proof}
Modulo $\ell$ one has $b\equiv\eta a$, and $\eta$ is a unit because
$\gcd(N\sigma,m)=1$.
\end{proof}

We next isolate the approximation estimate in a form that applies to every
vector of the Harvey--Hittmeir lattice, rather than only to a particular
short vector.

\begin{lemma}
\label{lem:norm-to-hh-error}
Let $m_0,\sigma_0>0$ and assume
\[
  \sigma_0\le p<\left(1+\frac1{m_0}\right)\sigma_0.
\]
Define
\begin{equation}
\label{eq:hh-linear-map}
  T(a,b)=(c,d)
  =\bigl(Na,\,-Nm_0a+m_0\sigma_0^2b\bigr).
\end{equation}
Then, for every $a,b\in\Z$,
\begin{equation*}
%\label{eq:norm-to-error}
  \left|
  aq+bp-\left(a\frac N{\sigma_0}+b\sigma_0\right)
  \right|
  \le
  \frac{\sqrt2\,\|T(a,b)\|}{m_0^2\sigma_0}.
\end{equation*}
Moreover, the image lattice $L'=T(L_\sigma)$ has determinant
\begin{equation}
\label{eq:transformed-det}
  D_L:=\det L'=Nmm_0\sigma_0^2.
\end{equation}
\end{lemma}

\begin{proof}
Write $p=\sigma_0(1+t_0)$, where $0\le t_0<1/m_0$.  As in the proof of
Harvey--Hittmeir Proposition~3.3, there is $0\le\delta<1$ such that
\[
  q=\frac N{\sigma_0}(1-t_0+\delta t_0^2).
\]
Consequently
\[
  aq+bp-\left(a\frac N{\sigma_0}+b\sigma_0\right)
  =t_0\left(-a\frac N{\sigma_0}+b\sigma_0\right)
   +\delta t_0^2a\frac N{\sigma_0}.
\]
For~\eqref{eq:hh-linear-map},
\[
  -a\frac N{\sigma_0}+b\sigma_0=\frac d{m_0\sigma_0},
  \qquad
  a\frac N{\sigma_0}=\frac c{\sigma_0}.
\]
Since $t_0<1/m_0$,
\[
  |\text{error}|
  \le\frac{|c|+|d|}{m_0^2\sigma_0}
  \le\frac{\sqrt2\,\|(c,d)\|}{m_0^2\sigma_0}.
\]
Applying $T$ to the basis $(1,\eta),(0,m)$ gives determinant
$Nmm_0\sigma_0^2$.
\end{proof}

The next lemma computes the local digit of the actual integer index used in
Harvey--Hittmeir Algorithm~4.3.  The intermediate quotients are
$m$-integral rationals, whereas the final index is an integer.

\begin{lemma}
\label{lem:quadratic-hh-index}
Assume that $m$ is squarefree, let $(a,b)\in L_\sigma$, and write
$p=\sigma+mk_p$ with $k_p\in\Z$.  Put
\[
  C=a\frac N\sigma+b\sigma,
  \qquad
  H=\frac{-aN/\sigma+b\sigma}{m}\in\Z_{(m)}.
\]
Let $\tau\in\{0,\ldots,m^2-1\}$ be the integer representative satisfying
$\tau\equiv C\pmod{m^2}$, and let $\Lambda\ge0$ be an integer.  Define
\[
  \tau_0=
  \left\lfloor
    \frac{aN/\sigma_0+b\sigma_0}{m^2}
  \right\rfloor
\]
and suppose that the Harvey--Hittmeir index $i$ is defined by
\begin{equation}
\label{eq:hh-index-definition}
  aq+bp=m^2(\tau_0-\Lambda+i)+\tau.
\end{equation}
Then for every prime $\ell\mid m$,
\begin{equation}
\label{eq:actual-hh-quadratic-index}
  i\equiv
  \frac{aN}{\sigma^3}k_p^2
  +\frac H\sigma k_p
  +R_\ell-\tau_0+\Lambda
  \pmod\ell,
\end{equation}
where
\[
  R_\ell\equiv\frac{C-\tau}{m^2}\pmod\ell
\]
is well defined in $\Z_{(m)}$ and is computable from
$C\pmod{m^2\ell}$.  If $\ell\nmid a$, the right-hand side of
\eqref{eq:actual-hh-quadratic-index} is a quadratic polynomial in
$k_p\pmod\ell$ with nonzero leading coefficient.
\end{lemma}

\begin{proof}
In $\Z_{(m)}$, put $u=mk_p/\sigma$.  Then
$p=\sigma(1+u)$ and $u\equiv0\pmod m$.  Modulo $m^3$,
\[
  q=\frac N\sigma(1+u)^{-1}
  \equiv\frac N\sigma(1-u+u^2).
\]
Hence
\[
  aq+bp-C
  \equiv
  \frac{mk_p}{\sigma}\left(-a\frac N\sigma+b\sigma\right)
  +\frac{m^2k_p^2aN}{\sigma^3}
  \pmod{m^3}.
\]
Since $-aN/\sigma+b\sigma=mH$, division by $m^2$ in $\Z_{(m)}$ gives
\begin{equation}
\label{eq:localized-next-digit}
  \frac{aq+bp-C}{m^2}
  \equiv
  \frac H\sigma k_p+\frac{aN}{\sigma^3}k_p^2
  \pmod m.
\end{equation}
By~\eqref{eq:hh-index-definition},
\[
  i=\frac{aq+bp-\tau}{m^2}-\tau_0+\Lambda
   =\frac{aq+bp-C}{m^2}+\frac{C-\tau}{m^2}-\tau_0+\Lambda.
\]
Reducing this identity modulo $\ell$ and using
\eqref{eq:localized-next-digit} proves
\eqref{eq:actual-hh-quadratic-index}.  To compute $R_\ell$, calculate
$C$ modulo $m^2\ell$, subtract the known integer $\tau$, and divide by
$m^2$.  The modular inverse of $\sigma$ exists because $\gcd(\sigma,m)=1$.
Finally, if $\ell\nmid a$, then $aN\sigma^{-3}$ is a unit modulo $\ell$.
\end{proof}

We construct a fixed forbidden set containing asymptotically one quarter of the local residues.

\begin{proposition}
\label{prop:quarter-forbidden-set}
Let $\ell$ be an odd prime, let $\chi$ be its quadratic character, and put
\[
  \mathcal R_\ell=\{x\in\F_\ell^*: \chi(x)=1\},
  \qquad
  \mathcal N_\ell=\{x\in\F_\ell: \chi(x)=-1\},
\]
\[
  F_\ell=\mathcal N_\ell\cap(1+\mathcal R_\ell).
\]
Then
\begin{equation}
\label{eq:quarter-size}
  |F_\ell|=\frac{\ell-\chi(-1)}4.
\end{equation}
For every quadratic polynomial
\[
  Q(X)=AX^2+BX+C\in\F_\ell[X],\qquad A\ne0,
\]
there is a computable translation $c(Q)\in\F_\ell$ such that
\begin{equation*}
%\label{eq:quarter-avoidance}
  (Q(\F_\ell)+c(Q))\cap F_\ell=\varnothing.
\end{equation*}
\end{proposition}

\begin{proof}
Expanding the indicator of $F_\ell$ gives
\[
  4|F_\ell|
  =\sum_{x\ne0,1}(1-\chi(x))(1+\chi(x-1)).
\]
The three character sums that occur are
\[
  \sum_{x\ne0,1}\chi(x)=-1,
  \qquad
  \sum_{x\ne0,1}\chi(x-1)=-\chi(-1),
\]
\[
  \sum_{x\in\F_\ell}\chi(x(x-1))=-1.
\]
Substitution yields~\eqref{eq:quarter-size}.

Completing the square gives
\[
  Q(\F_\ell)=u+A(\{0\}\cup\mathcal R_\ell)
\]
for a computable $u\in\F_\ell$.  If $\chi(A)=1$, take $c(Q)=-u$.
Then the translated image is $\{0\}\cup\mathcal R_\ell$ and its complement
contains $\mathcal N_\ell$, hence $F_\ell$.  If $\chi(A)=-1$, take
$c(Q)=1-u$.  Then the translated image is
$1+(\{0\}\cup\mathcal N_\ell)$ and its complement contains
$1+\mathcal R_\ell$, again containing $F_\ell$.
\end{proof}

\begin{proposition}
\label{prop:quarter-optimality}
Let $d\in\F_\ell^*$.  Then
\begin{equation}
\label{eq:quadratic-class-intersection}
  |\mathcal N_\ell\cap(d+\mathcal R_\ell)|
  =\frac{\ell-1+\chi(d)-\chi(-d)}4.
\end{equation}
Consequently the largest set that can be forced into the complement of both
quadratic-image square classes by translations alone has cardinality
$(\ell-1)/4$ when $\ell\equiv1\pmod4$ and $(\ell+1)/4$ when
$\ell\equiv3\pmod4$.  Proposition~\textup{\ref{prop:quarter-forbidden-set}} attains
this bound.
\end{proposition}

\begin{proof}
Suppose that a fixed set $F$ can be placed in the complement of both
quadratic-image square classes by translations.  After normalizing the square
image to $\{0\}\cup\mathcal R_\ell$, its complement is
$\mathcal N_\ell$.  After normalizing the nonsquare image to
$\{0\}\cup\mathcal N_\ell$, its complement is
$\mathcal R_\ell$.  Thus for suitable translations $c_1,c_2$ one must
have
\[
  F\subseteq(c_1+\mathcal N_\ell)\cap(c_2+\mathcal R_\ell).
\]
Translating by $-c_1$ reduces the problem to bounding
$|\mathcal N_\ell\cap(d+\mathcal R_\ell)|$, where
$d=c_2-c_1$.  The case $d=0$ gives the empty intersection, so assume
$d\ne0$.  Then
\[
  4|\mathcal N_\ell\cap(d+\mathcal R_\ell)|
  =\sum_{x\ne0,d}(1-\chi(x))(1+\chi(x-d)).
\]
Using
\[
  \sum_{x\ne0,d}\chi(x)=-\chi(d),
  \quad
  \sum_{x\ne0,d}\chi(x-d)=-\chi(-d),
\]
and
\[
  \sum_x\chi(x(x-d))=-1
\]
gives~\eqref{eq:quadratic-class-intersection}.  Maximizing over $d$ yields
the stated values, which agree with~\eqref{eq:quarter-size}.
\end{proof}

The next elementary observation is what allows the same set of selected
primes to be usable for every regular giantstep.

\begin{lemma}
\label{lem:upper-prime-avoidance}
Let $w_1=(a_1,b_1)$ and $w_2=(a_2,b_2)$ be a basis of $L_\sigma$.  Then
$\gcd(a_1,a_2)=1$.  Let $\mathcal S$ be a set of $t$ prime divisors of $m$
such that every $\ell\in\mathcal S$ satisfies $\ell>2t$, where $t\ge1$.  Then there exists
$k\in\{0,\ldots,2t-1\}$ such that
\begin{equation}
\label{eq:all-selected-good}
  \ell\nmid a_1+ka_2
  \qquad(\ell\in\mathcal S).
\end{equation}
For the same $k$, no $\ell\in\mathcal S$ divides the second coordinate
$b_1+kb_2$.
\end{lemma}

\begin{proof}
The projection $(a,b)\mapsto a$ maps $L_\sigma$ onto $\Z$, because
$(1,\eta)\in L_\sigma$.  Since $w_1,w_2$ is a basis, this image is
$a_1\Z+a_2\Z$, proving $\gcd(a_1,a_2)=1$.

Fix $\ell\in\mathcal S$.  If $\ell\nmid a_2$, exactly one residue class
$k\pmod\ell$ makes $a_1+ka_2$ divisible by $\ell$.  Because
$0\le k<2t<\ell$, this excludes at most one actual value of $k$.  If
$\ell\mid a_2$, then $\ell\nmid a_1$ and no $k$ is excluded.  Thus the $t$
primes exclude at most $t$ of the $2t$ possible values, so an admissible $k$
exists.  The assertion for the second coordinate follows from
Lemma~\ref{lem:local-divisibility-dichotomy}.
\end{proof}

We now combine the preceding lemma with Gauss reduction.  The parameter $Y$
will later be chosen to balance the regular and skinny branches.

\begin{proposition}
\label{prop:quarter-gauss-dichotomy}
Let $z_1,z_2$ be a Gauss-reduced basis of $L'=T(L_\sigma)$, ordered by
$\|z_1\|\le\|z_2\|$, and let $w_1,w_2$ be the inverse images in $L_\sigma$.
Let $\mathcal S$ satisfy the hypotheses of
Lemma~\textup{\ref{lem:upper-prime-avoidance}}, and let $Y\ge1$.  Then one of the
following holds.
\begin{enumerate}[label=\textup{(\roman*)}]
\item If $\|z_2\|\le Y\sqrt{D_L}$, there exists $0\le k<2t$ for which
$w=w_1+kw_2=(a,b)$ satisfies~\eqref{eq:all-selected-good} and
\begin{equation*}
%\label{eq:regular-vector-bound}
  \|T(w)\|\le(1+2t)Y\sqrt{D_L}.
\end{equation*}
\item If $\|z_2\|>Y\sqrt{D_L}$, then the vector $w_1=(a,b)$ satisfies
\begin{equation*}
%\label{eq:skinny-vector-bound}
  \|T(w_1)\|<\frac{2}{\sqrt3\,Y}\sqrt{D_L}.
\end{equation*}
\end{enumerate}
\end{proposition}

\begin{proof}
For a Gauss-reduced basis the angle between $z_1$ and $z_2$ lies in
$[\pi/3,2\pi/3]$.  Hence
\[
  D_L\ge\frac{\sqrt3}{2}\|z_1\|\|z_2\|.
\]
In the first case apply Lemma~\ref{lem:upper-prime-avoidance} and use
\[
  \|z_1+kz_2\|
  \le\|z_1\|+2t\|z_2\|
  \le(1+2t)Y\sqrt{D_L}.
\]
In the second case divide the product bound by
$\|z_2\|>Y\sqrt{D_L}$.
\end{proof}

For later use, define
\[
  M=\prod_{\ell\in\mathcal S}\ell,
  \qquad
  \rho=\prod_{\ell\in\mathcal S}
  \left(1-\frac{|F_\ell|}{\ell}\right).
\]
For each regular giantstep, Lemma~\ref{lem:quadratic-hh-index} and
Proposition~\ref{prop:quarter-forbidden-set} produce local translations
$c_\ell$, and CRT produces $c\pmod M$.  After replacing the giantstep $v$ by
$\alpha^cv$, where $\alpha=\beta^{m^2}$, the true exponent $i+c$ belongs to
the single candidate-independent set
\begin{equation}
\label{eq:universal-quarter-set}
  E_R(L)=\{0\le r\le L:r\bmod\ell\notin F_\ell
                    \text{ for every }\ell\in\mathcal S\}.
\end{equation}
The following elementary estimate will be used repeatedly.
\begin{equation}
\label{eq:finite-interval-quarter-count}
  |E_R(L)|\le\rho(L+M).
\end{equation}
Indeed, every complete block of length $M$ contains exactly $\rho M$
retained residues, and there is at most one incomplete block.

We require arbitrary-point multipoint evaluation over the composite ring
$\Z/N\Z$.  We record a self-contained form sufficient for the algorithm.

\begin{lemma}
\label{lem:ring-multipoint}
Let $f\in(\Z/N\Z)[X]$ have degree at most $n$, and let
$x_1,\ldots,x_s\in\Z/N\Z$.  All values $f(x_j)$ can be computed
deterministically in
\begin{equation*}
%\label{eq:ring-multipoint-cost}
  O((n+s)\log^{3+o(1)} N)
\end{equation*}
bit operations whenever $n+s=N^{O(1)}$.
\end{lemma}

\begin{proof}
Construct the usual subproduct tree from the monic linear polynomials
$X-x_j$.  At every level the sum of the degrees is $O(s)$, so using fast
integer multiplication and Kronecker substitution, all products on one level
cost $O(s\log^{2+o(1)}N)$ bit operations.  There are $O(\log s)$ levels.

Before descending, compute the root remainder
$F_{\rm root}=f\bmod P_{\rm root}$.  This is division by a monic polynomial
of degree $s$ and costs
$O((n+s)\log^{2+o(1)}N)$ bit operations by the same fast multiplication
model.  To descend the remainder tree one again divides only by monic
polynomials.  Such division is valid over any commutative ring.  In the fast
reversal method the reversal of a monic divisor has constant term $1$, and
hence has a truncated formal-power-series inverse obtained by Newton
iteration using only additions and multiplications in $\Z/N\Z$.  Thus no
inversion of a nonunit modulo $N$ is required.  The same level-by-level
complexity bound applies to the remainder tree.  Since $n+s=N^{O(1)}$,
$\log(n+s)=O(\log N)$, and the stated $O((n+s)\log^{3+o(1)}N)$ bound follows.
\end{proof}

\begin{lemma}
\label{lem:sparse-collision}
Let $N=pq$ with distinct primes $p,q$, let $\alpha\in(\Z/N\Z)^*$, let $E$ be a
finite set of nonnegative integers bounded by $N^{O(1)}$, and let
$v_1,\ldots,v_n\in\Z/N\Z$.  Suppose exact equalities
$v_h=\alpha^e\pmod N$ with $e\in E$ have first been removed and treated
separately.  Assume $n+|E|=N^{O(1)}$.  If a collision modulo $p$ or $q$ exists between the
sets $\{\alpha^e:e\in E\}$ and $\{v_h\}$, a proper factor of $N$ can be found in
\begin{equation}
\label{eq:sparse-collision-cost}
  O((n+|E|)\log^{3+o(1)}N)
\end{equation}
bit operations.
\end{lemma}

\begin{proof}
Form
\[
  f(X)=\prod_{h=1}^n(X-v_h)
\]
and evaluate $f(\alpha^e)$ for all $e\in E$ using
Lemma~\ref{lem:ring-multipoint}.  If
$1<\gcd(f(\alpha^e),N)<N$, a factor has been found.  If the gcd equals $N$,
then, for a semiprime $N=pq$, at least one factor
$\alpha^e-v_h$ is divisible by $p$ and at least one is divisible by $q$.
Because exact modulo-$N$ matches were removed, these cannot be the same
factor in the only case where every gcd would be trivial.  Scanning
$\gcd(\alpha^e-v_h,N)$ therefore yields a proper factor.  The algorithm
terminates as soon as this case occurs, so such a full scan is performed at
most once.  Product construction, multipoint evaluation, and the possible
final scan satisfy the stated bound.
\end{proof}

The next lemma replaces the linear factor-order certification loop used in
Harvey--Hittmeir Algorithm~4.3 by a square-root-cost deterministic test.

\begin{lemma}
\label{lem:factor-order-certification}
Let $N=pq$ with distinct primes $p,q$, let $\alpha\in(\Z/N\Z)^*$, and let the integer $L\ge1$ satisfy
$\operatorname{ord}_N(\alpha)>L$.  There is a deterministic algorithm costing
$L^{1/2}(\log N)^{O(1)}$ bit operations that either finds a nontrivial factor
of $N$ or certifies
\[
  \operatorname{ord}_p(\alpha)>L,
  \qquad
  \operatorname{ord}_q(\alpha)>L.
\]
\end{lemma}

\begin{proof}
Put $r=\lceil\sqrt L\rceil$ and $s=\lfloor L/r\rfloor$.  Compare
\[
  \alpha^i\quad(0\le i<r)
\]
with
\[
  \alpha^{jr}\quad(1\le j\le s)
\]
using a product-polynomial collision test.  Every exponent difference
$jr-i$ lies in $[1,L]$, so a collision modulo $N$ is impossible under the
hypothesis $\operatorname{ord}_N(\alpha)>L$.  If
$d=\operatorname{ord}_p(\alpha)\le sr$, write
$j=\lceil d/r\rceil$ and $i=jr-d$, which produces a collision modulo $p$.
Use the same product-polynomial/GCD collision procedure as in
Lemma~\ref{lem:sparse-collision}.  A gcd equal to $N$ can only mean that
one factor collides modulo $p$ and another modulo $q$.  Because an exact
modulo-$N$ collision is impossible when the exponent difference is at most
$L<\operatorname{ord}_N(\alpha)$, scanning the individual differences then
yields a proper factor.  The same argument applies with $p$ and $q$
interchanged.

There remain fewer than $r$ possible orders $d$ with $sr<d\le L$.  Testing
$\gcd(\alpha^d-1,N)$ for these values completes the certification.  If one
of these gcds were $N$, then both factor orders would divide $d\le L$, which
would contradict $\operatorname{ord}_N(\alpha)>L$, so any nontrivial tail
gcd is proper.  The two lists and the tail contain $O(\sqrt L)$ entries, and
the product-polynomial operations have total cost $O(\sqrt L\log^{3+o(1)}N)$.  Computing the tail powers by repeated squaring costs $O(\sqrt L\log L\log^{1+o(1)}N)$.  Since $L<\operatorname{ord}_N(\alpha)<N$, both bounds have the asserted form.
\end{proof}

We can now state the modified search bound.  It is an analogue of
Harvey--Hittmeir Proposition~4.4 with the final collision term multiplied by
the quarter-sieve factor arising from the regular/skinny dichotomy.

\begin{proposition}
\label{prop:modified-hh-search}
Let $N=pq$ with distinct primes $p<q$, let $m$ be squarefree with $\gcd(m,N)=1$, and
let $m_0\ge1$ be an integer.  Assume $m,m_0\le N$.  Let $\mathcal S$ be a set of $t$ prime divisors of $m$ such
that every $\ell\in\mathcal S$ satisfies $\ell>2t$, and define $M$ and
$\rho$ as above.  Put
\[
  X=mm_0,
  \qquad
  \lambda_0=\frac{N^{1/2}}{X^{3/2}},
  \qquad
  R=(t\rho)^{1/2},
  \qquad
  Y=R^{-1}.
\]
Assume $t\ge1$, $t\rho\le1$, and $\lambda_0\ge1$, so $Y\ge1$.  Assume that an element $\beta\in(\Z/N\Z)^*$ has been found such that, for
$\alpha=\beta^{m^2}$,
\[
  \operatorname{ord}_N(\alpha)>L_*,
\]
where
\begin{equation}
\label{eq:max-exponent-span}
  L_*=\left\lceil C\lambda_0\sqrt{t/\rho}\right\rceil+M
\end{equation}
for a sufficiently large absolute constant $C$.  Then a deterministic
modification of Harvey--Hittmeir Algorithm~$4.3$ factors $N$, and it runs in
\begin{equation}
\label{eq:modified-search-bound}
  O\left(
    m(\log N)^{O(1)}
    +\varphi(m)m_0\log^4N
    +R\frac{N^{1/2}}{(mm_0)^{3/2}}\log^{3+o(1)}N
    +\rho M\log^{3+o(1)}N
  \right)
\end{equation}
bit operations, apart from the one-time factor-order certification of
Lemma~\textup{\ref{lem:factor-order-certification}}.
\end{proposition}

\begin{proof}
The residue loop visits every $\sigma\in\{1,\ldots,m\}$ with $\gcd(\sigma,m)=1$.  For each residue, start at $\sigma_0=1$ and repeatedly replace $\sigma_0$ by $\lceil(1+1/m_0)\sigma_0\rceil$ while $\sigma_0<\sqrt N$.  These are the loops of Harvey--Hittmeir Algorithm~4.3.  There are
$O(\varphi(m)m_0\log N)$ pairs $(\sigma,\sigma_0)$.
For each pair compute a Gauss-reduced basis of $L'=T(L_\sigma)$.

If the second reduced vector has norm at most $Y\sqrt{D_L}$, use
Proposition~\ref{prop:quarter-gauss-dichotomy}(i) to choose
$w=w_1+kw_2$ with $k<2t$, good at every prime in $\mathcal S$.  By
Lemmas~\ref{lem:norm-to-hh-error} and the elementary inequality
$|\lfloor x\rfloor-\lfloor y\rfloor|\le|x-y|+1$, the correct integer
$\lfloor(aq+bp)/m^2\rfloor$ lies in an interval of radius
\[
  \Lambda_R\ll tY\lambda_0+1.
\]
Use this common upper bound in the definition of the Harvey--Hittmeir index.
Lemma~\ref{lem:quadratic-hh-index} and
Proposition~\ref{prop:quarter-forbidden-set} give, for each
$\ell\in\mathcal S$, a computable translation $c_\ell$.  CRT gives
$c\pmod M$.  Multiplying the giantstep by $\alpha^c$ moves its true exponent
from $i$ to $i+c$, and all true shifted exponents belong to the single set
$E_R(2\Lambda_R+M)$ from~\eqref{eq:universal-quarter-set}.  By
\eqref{eq:finite-interval-quarter-count},
\[
  |E_R(2\Lambda_R+M)|
  \ll \rho tY\lambda_0+\rho M+1
  =R\lambda_0+\rho M+1.
\]
All regular giantsteps can therefore be processed by one invocation of
Lemma~\ref{lem:sparse-collision}.

If the second reduced vector has norm greater than $Y\sqrt{D_L}$, use
Proposition~\ref{prop:quarter-gauss-dichotomy}(ii).  The corresponding radius
satisfies
\[
  \Lambda_S\ll Y^{-1}\lambda_0+1=R\lambda_0+1.
\]
These skinny giantsteps are processed by the ordinary dense collision search
on $0\le i\le2\Lambda_S$.

Before either collision search, invoke
Lemma~\ref{lem:factor-order-certification} with the common exponent span
$L_*$.  Hence
\[
  \operatorname{ord}_p(\alpha)>L_*,
  \qquad
  \operatorname{ord}_q(\alpha)>L_*.
\]
We now verify correctness directly for the unique pair
$(\bar\sigma,\bar\sigma_0)$ satisfying
$p\equiv\bar\sigma\pmod m$ and
$\bar\sigma_0\le p<(1+1/m_0)\bar\sigma_0$, whose existence and uniqueness follow from the loops just described.  Indeed, for successive values $s$ and $s^{\prime}=\lceil(1+1/m_0)s\rceil$, the condition $s\le p<s^{\prime}$ is equivalent to $s\le p<(1+1/m_0)s$ because $p$ is an integer.  These intervals partition the integers below $\sqrt N$, and $p<\sqrt N$.

Let $(a,b)$ be the vector chosen by the appropriate branch, let $\Lambda$
be the corresponding common radius, and put
\[
  j=m^2(\tau_0-\Lambda)+\tau.
\]
The approximation bound proved above gives an integer
$0\le i\le2\Lambda$ such that
\[
  aq+bp=m^2i+j.
\]
Since $N=pq$, reduction modulo $p-1$ gives
\[
  aq+bp\equiv aN+b\pmod{p-1}.
\]
Therefore Fermat's little theorem implies
\[
  \alpha^i=\beta^{m^2i}
  \equiv \beta^{aN+b-j}\pmod p.
\]
In the skinny branch this is exactly a collision between the dense babystep
set and the corresponding giantstep.  In the regular branch, if $c$ is the
CRT translation attached to this candidate and
$v'=\alpha^c\beta^{aN+b-j}$, then
\[
  \alpha^{i+c}\equiv v'\pmod p.
\]
By Lemma~\ref{lem:quadratic-hh-index} and
Proposition~\ref{prop:quarter-forbidden-set}, the exponent $i+c$ belongs to
the universal sparse set $E_R(2\Lambda_R+M)$.  Thus the correct candidate
produces a collision modulo $p$ in one of the two searches.

It remains to justify removal of exact modulo-$N$ matches.  Consider the
correct giantstep and let $e_*$ denote its true exponent ($e_*=i$ in the
skinny branch and $e_*=i+c$ in the regular branch).  If the same giantstep
were exactly equal modulo $N$ to $\alpha^e$ for another tested exponent
$e$, then modulo $p$ one would have $\alpha^e=\alpha^{e_*}$.  Both exponents
lie in $[0,L_*]$, so the certified inequality
$\operatorname{ord}_p(\alpha)>L_*$ forces $e=e_*$.  Hence an exact match
for the correct giantstep is necessarily the true match.  Apply the recovery test of~\cite[Lemma~2.4]{HarveyHittmeir2022} to
\[
  u=m^2i+j=aq+bp
\]
to recover the factors.  To make the recovery valid for arbitrary signs, if $ab\ne0$ compute the integer roots of $Z^2-uZ+abN$.  For the true candidate they are $aq$ and $bp$.  Test both assignments, divide by $a$ and $b$, and accept only positive integers greater than $1$ whose product is $N$.  This uses exact integer square roots and divisions on $O(\log N)$-bit integers.  No nonvanishing assumption on both coefficients is needed here.  If $a=0$ and $b\ne0$, then $p=u/b$, while if $b=0$ and
$a\ne0$, then $q=u/a$.  The case $a=b=0$ cannot occur because the chosen
lattice vector is nonzero.  In the regular branch, the selected-prime
avoidance already forces $a\ne0$.  Consequently removing any other exact matches cannot
discard the correct collision.  If no exact match recovers the factors, the
correct giantstep remains in the appropriate collision list, and
Lemma~\ref{lem:sparse-collision} (or the ordinary dense collision routine in
the skinny branch) finds a proper factor from its collision modulo $p$.

It remains to account for constructing the sparse evaluation set and its
powers.  Let $A_M\subseteq\{0,\ldots,M-1\}$ be the retained residue
classes modulo $M$.  We may list $A_M$ by scanning one complete residue block, which costs $M(\log N)^{O(1)}$ bit operations.  For each $r\in A_M$ with $r\le L$, output $r+jM$ for $0\le j\le\lfloor(L-r)/M\rfloor$.  This enumerates $E_R(L)$ in $O(M+|E_R(L)|)$ arithmetic operations, without scanning the excluded exponents.  Compute each retained power by repeated squaring.  Since $L\le L_*=N^{O(1)}$ under the stated assumptions, this costs $O(|E_R(L)|\log^{2+o(1)}N)$ bit operations.  As $M\le m$, the enumeration and exponentiation costs are absorbed by the first term and the collision terms in~\eqref{eq:modified-search-bound}.

The per-candidate work consists of two-dimensional lattice reduction, marking
at most $t=O(\log N)$ forbidden values of $k$, computing the local translations,
one CRT reconstruction, and modular exponentiation.  All operands have
$O(\log N)$ bits, so this is $O(\log^3N)$ bit operations per candidate and
contributes $O(\varphi(m)m_0\log^4N)$ in total.  The outer residue loop has
the same $m(\log N)^{O(1)}$ cost as in Harvey--Hittmeir.  The regular and
skinny collision-set sizes were bounded above by
$O(R\lambda_0+\rho M+1)$, and
Lemma~\ref{lem:sparse-collision}, together with the preceding enumeration bound, gives the final two terms in
\eqref{eq:modified-search-bound}.
\end{proof}

We finally optimize the auxiliary modulus.  We use the unrestricted
large-order theorem of Harvey and Hittmeir~\cite{HarveyHittmeir2026}, which finds, for every $1\le D<N-1$, either a nontrivial divisor of $N$ or an element of order exceeding $D$ in $O(D^{1/2}\log D\log N/\sqrt{\log\log D})$ bit operations.

\begin{theorem}
\label{thm:factorization-main}
There is a deterministic algorithm that factors every positive integer $N$
in
\begin{equation}
\label{eq:main-factorization-bound}
  N^{1/5}
  \exp\left(
    -\left(\frac{\log(4/3)}{50}+o(1)\right)
     \frac{\log N}{\log\log N}
  \right)
\end{equation}
bit operations.
\end{theorem}

Theorem~\ref{thm:factorization-main} improves the general Harvey--Hittmeir factorization bound with logarithmic exponent $16/5$~\cite{HarveyHittmeir2022}.  For balanced semiprimes, its secondary factor is also asymptotically smaller than that in the rank-$3$ bound with logarithmic exponent $13/5$~\cite{GaoFengHuPan2025}.  The unrestricted large-order theorem of Harvey and Hittmeir~\cite{HarveyHittmeir2026} is used as an ingredient in the proof.  The theorem applies to general inputs after the standard preliminary reduction.

\begin{proof}
First remove powers of $2$, test for a perfect square, and apply deterministic primality testing.  A square is handled by factoring its integer square root.  For the remaining composite inputs it is enough to treat products of two distinct primes because the standard preliminary reduction used in the
exponent-one-fifth algorithms~\cite{Harvey2021,HarveyHittmeir2022} factors
an integer with at least three prime factors, counted with multiplicity, in
$N^{1/6+o(1)}$ time, which is negligible compared with
\eqref{eq:main-factorization-bound}, and recursion contributes only polynomial
factors in $\log N$.  Thus for the main argument we may write $N=pq$ with distinct primes $p<q$.

Put $L=\log N$, $L_2=\log\log N$, $E=L/L_2$, and $a=\log(4/3)$.  Choose a fixed constant
\begin{equation}
\label{eq:s-choice}
  s>\frac{3a}{50}
\end{equation}
and set
\begin{equation*}
%\label{eq:primorial-endpoint}
  B=\frac{L}{10}-sE.
\end{equation*}
Let $m=\prod_{\ell\le B}\ell$.  We use the classical quantitative prime number theorem in the form
\[
  \vartheta(x)=x+O\!\left(xe^{-c_1\sqrt{\log x}}\right),
  \qquad
  \pi(x)=\operatorname{Li}(x)
   +O\!\left(xe^{-c_1\sqrt{\log x}}\right)
\]
for some absolute $c_1>0$, as stated in, for example,
Davenport~\cite{Davenport2000}.  Since $B\asymp L$, the error term is
$o(B/\log B)=o(E)$.  Hence
\begin{equation}
\label{eq:m-endpoint}
  \log m=\vartheta(B)=\frac{L}{10}-sE+o(E).
\end{equation}
Before proceeding compute $g=\gcd(N,m)$.  If $1<g<N$, we have already
found a nontrivial factor.  Henceforth assume $\gcd(N,m)=1$.

Let $y=(2+\varepsilon)B/\log B$ for any fixed $\varepsilon>0$, and let $\mathcal S$ be the primes in
$(y,B]$.  The same quantitative prime number theorem, together with the
standard expansion $\operatorname{Li}(B)=B/\log B+O(B/\log^2B)$, gives
\begin{equation*}
%\label{eq:t-endpoint}
  t=|\mathcal S|=\frac{B}{\log B}+o(E)
   =\frac1{10}E+o(E).
\end{equation*}
Moreover $2t<y$ for sufficiently large $N$, so
Lemma~\ref{lem:upper-prime-avoidance} applies.

By Proposition~\ref{prop:quarter-forbidden-set},
\[
  \rho
  =\prod_{\ell\in\mathcal S}
   \left(\frac34+O(1/\ell)\right)
  =\exp\left(-\frac a{10}E+o(E)\right).
\]
Consequently
\begin{equation*}
%\label{eq:R-endpoint}
  R=(t\rho)^{1/2}
  =\exp\left(-\frac a{20}E+o(E)\right).
\end{equation*}

Let $X_0$ denote the positive real number obtained by balancing the second
and third terms of Proposition~\ref{prop:modified-hh-search}, including the
displayed powers of $\log N$, and choose
\[
  m_0=\left\lceil\frac{X_0}{m}\right\rceil,
  \qquad X=mm_0.
\]
The estimates below show $X_0/m=N^{1/10+o(1)}$, so rounding changes $X$ by a
relative factor $1+o(1)$ and has no effect on any displayed secondary
asymptotic.  Apart from ordinary powers of $\log N$, the two dominant terms
are
\[
  \frac{\varphi(m)}mX
  \quad\text{and}\quad
  RN^{1/2}X^{-3/2}.
\]
Mertens' theorem gives
$\varphi(m)/m=(\log\log N)^{-1+o(1)}$.  Balancing the displayed terms,
including the powers of $\log N$ in
\eqref{eq:modified-search-bound}, yields
\begin{equation}
\label{eq:X-endpoint}
  \log X
  =\frac L5-\frac a{50}E+o(E).
\end{equation}
The first two terms of~\eqref{eq:modified-search-bound} and the balanced
collision term therefore contribute
\begin{equation}
\label{eq:search-endpoint}
  N^{1/5}\exp\left(-\frac a{50}E+o(E)\right).
\end{equation}
The direct $m$ term is only $N^{1/10+o(1)}$.

The selected-prime modulus satisfies
$\log M=B-y+o(B)=L/10+o(L)$, so
$M=N^{1/10+o(1)}$.  Hence the boundary term $\rho M\log^{3+o(1)}N$ is also
polynomially smaller than~\eqref{eq:search-endpoint}.

The estimates above give $t\ge1$, $t\rho<1$, $\lambda_0\ge1$, and $m,m_0\le N$ for sufficiently large $N$, as required by Proposition~\ref{prop:modified-hh-search}.  It remains to justify its order assumption.  From
\eqref{eq:X-endpoint},
\[
  \lambda_0=N^{1/2}X^{-3/2}
  =N^{1/5}
    \exp\left(\frac{3a}{100}E+o(E)\right).
\]
Since $\sqrt{t/\rho}=\exp(aE/20+o(E))$, the maximum exponent span in
\eqref{eq:max-exponent-span} satisfies
\begin{equation}
\label{eq:Lstar-endpoint}
  L_*
  =N^{1/5}\exp\left(\frac{2a}{25}E+o(E)\right).
\end{equation}
Ask the large-order algorithm of
Harvey--Hittmeir~\cite[Theorem~1.1]{HarveyHittmeir2026} for an element
$\beta$ of order exceeding
\[
  D=m^2L_*.
\]
By~\eqref{eq:m-endpoint} and~\eqref{eq:Lstar-endpoint},
\[
  \log D
  =\frac{2L}{5}
   +\left(-2s+\frac{2a}{25}\right)E+o(E),
\]
so $D<N-1$ for all sufficiently large $N$.  The cost of finding $\beta$ is
\[
  N^{1/5}
  \exp\left(\left(-s+\frac a{25}\right)E+o(E)\right).
\]
By~\eqref{eq:s-choice}, this is
$o\bigl(N^{1/5}\exp(-aE/50)\bigr)$.  Moreover,
\[
  \operatorname{ord}_N(\beta^{m^2})
  =\frac{\operatorname{ord}_N(\beta)}
         {\gcd(\operatorname{ord}_N(\beta),m^2)}
  >\frac D{m^2}=L_*.
\]
Lemma~\ref{lem:factor-order-certification} then either factors $N$ or
certifies the corresponding order lower bounds modulo both prime factors, at
cost $L_*^{1/2}(\log N)^{O(1)}=N^{1/10+o(1)}$.

All remaining preprocessing and recursion costs are smaller than
\eqref{eq:search-endpoint}.  Since $a=\log(4/3)$, the bound
\eqref{eq:main-factorization-bound} follows.
\end{proof}

\begin{remark}
\label{rem:square-root-loss}
The constant $\log(4/3)/50$ is not claimed to be optimal.  By Proposition~\ref{prop:quarter-optimality}, the local quarter-forbidden set is optimal among universal forbidden sets obtained by translations and required to work for both quadratic square classes.  The factor $R=\sqrt{t\rho}$, whose logarithm differs from $(\log\rho)/2$ by $(\log t)/2$, comes from balancing the regular and skinny cases in Proposition~\ref{prop:quarter-gauss-dichotomy}.  A larger saving constant would therefore require either a better treatment of the skinny case or a different global construction.
\end{remark}

\section{Proof-of-concept computations and finite-size crossover}
\label{sec:computations}

We implemented a matched proof-of-concept comparison between the dense
Harvey--Hittmeir search and the quadratic-image-sieved search.  Both versions
use the same semiprime, auxiliary modulus, residue-class candidate grid,
group element, recovery test, and collision backend.  The implementation is
intended to check the structure of the argument and to compare the sizes of
the search sets.  It is not a tuned factorization program.

For every regular candidate, the code verifies the quadratic congruence for
the actual Harvey--Hittmeir babystep index, the translated membership
condition in the universal retained set, and the final collision congruence.
In the computations reported here all such checks passed.  The observed
retained densities also agree with the exact product
\[
 \rho=\prod_{\ell\in\mathcal S}
 \left(1-\frac{|F_\ell|}{\ell}\right)
\]
to the displayed precision.

To separate the effect of the primorial cutoff from the input size, we first
performed paired sweeps in which the same ten semiprimes were used for every
value of $B$.  Table~\ref{tab:B-sweep} reports the mean ratio of the number of
QIS babysteps to Harvey--Hittmeir babysteps, the ratio of collision-search
times, and the ratio of total proof-of-concept times.  Here ``total'' means
candidate generation plus collision search.  The absolute timings depend on
the implementation and hardware.  The ratios are included only to illustrate
the finite-size tradeoff.

\begin{table}[ht]
\centering
\caption{Paired proof-of-concept $B$-sweeps.  Each row contains ten
semiprimes.  Every factorization succeeded and the number of failed internal
QIS checks was zero.}
\label{tab:B-sweep}
\begin{tabular}{r r c c c c}
\hline
bits & $B$ & $\rho$ & baby ratio & search ratio & total ratio\\
\hline
72 & 7  & 0.5714 & 4.45 & 4.05 & 3.33\\
72 & 11 & 0.5195 & 4.08 & 4.52 & 3.62\\
72 & 13 & 0.3996 & 4.37 & 4.32 & 3.49\\
72 & 17 & 0.4278 & 5.18 & 0.76 & 1.21\\
80 & 7  & 0.5714 & 4.45 & 3.96 & 3.40\\
80 & 11 & 0.5195 & 4.08 & 3.28 & 2.86\\
80 & 13 & 0.3996 & 4.37 & 5.03 & 4.18\\
80 & 17 & 0.4278 & 4.67 & 1.94 & 1.42\\
\hline
\end{tabular}
\end{table}

These small inputs are far below the asymptotic range of
Theorem~\ref{thm:factorization-main}.  In particular, for fixed small $B$ the
local density gain is constant while the enlargement coming from the nearby
regular vector is also constant, so no asymptotic improvement should be
expected.  At the other extreme, taking $B$ too large causes the
residue-class candidate grid to dominate.  The 72-bit row with $B=17$
illustrates this directly because the collision-search phase is already smaller
than in the dense search, but the total proof-of-concept running time remains
larger because candidate generation is more expensive.

We next record where the asymptotic balance should become visible.  Put
\[
 E=\frac{\log N}{\log\log N},\qquad a=\log(4/3).
\]
The proof gives
\[
 t\sim \frac{E}{10},\qquad
 \rho=\exp\left(-\frac{a}{10}E+o(E)\right),
 \qquad
 R=(t\rho)^{1/2}.
\]
Thus the continuous approximation is
\begin{equation}
\label{eq:finite-crossover-proxy}
 R_{\rm cont}(E)
 =
 \sqrt{\frac{E}{10}}\,
 \exp\left(-\frac{a}{20}E\right).
\end{equation}
The larger solution of $R_{\rm cont}(E)=1$ is
$E=65.139\ldots$.  Solving
$\log N/\log\log N=65.139\ldots$ gives
\[
 \log_2 N=560.28\ldots.
\]
This number is only a continuous lower-order proxy.  At such sizes the prime
window used in the proof contains only a few primes, so the discreteness of
the selected-prime set is significant.

For a more realistic finite-size projection, choose for illustration
$s=0.018>3a/50$ in the endpoint condition and $\varepsilon=0.1$ in the
selected-prime window, and put
\[
 B_N=\frac{\log N}{10}-sE,\qquad
 y_N=(2+\varepsilon)\frac{B_N}{\log B_N}.
\]
Let $\mathcal S_N$ be the primes in $(y_N,B_N]$, let $t_N=|\mathcal S_N|$,
and form the exact finite product $\rho_N$ from these primes.  The resulting
values of $R_N=\sqrt{t_N\rho_N}$ are shown in
Table~\ref{tab:finite-crossover}.

\begin{table}[ht]
\centering
\caption{Finite-size projection for one admissible endpoint choice.  The
quantity $R_N=\sqrt{t_N\rho_N}$ is the regular/skinny reduction factor
appearing in the proof.}
\label{tab:finite-crossover}
\begin{tabular}{r r r c c}
\hline
bits & $B_N$ & $t_N$ & $\rho_N$ & $R_N$\\
\hline
256  & 17.1  & 2  & 0.588   & 1.085\\
384  & 25.8  & 3  & 0.416   & 1.118\\
512  & 34.4  & 3  & 0.416   & 1.117\\
560  & 37.6  & 4  & 0.315   & 1.122\\
640  & 43.1  & 5  & 0.240   & 1.095\\
700  & 47.1  & 6  & 0.178   & 1.035\\
768  & 51.7  & 6  & 0.178   & 1.035\\
800  & 53.9  & 7  & 0.135   & 0.971\\
1024 & 69.0  & 8  & 0.100   & 0.896\\
1536 & 103.7 & 13 & 0.0236  & 0.554\\
2048 & 138.4 & 17 & 0.00747 & 0.356\\
\hline
\end{tabular}
\end{table}

For this concrete finite model, $R_N$ first falls below $1$ between roughly
700 and 800 bits.  The $700$--$800$ bit range is therefore where this particular proxy changes from exceeding $1$ to being less than $1$.  Its values decrease further at $1024$ and $1536$--$2048$ bits.  These values should not be interpreted as measured
running-time crossovers.  Our present proof-of-concept implementation cannot
reach that range efficiently, and lower-order costs are not represented by
$R_N$.  They quantify the finite products associated with the asymptotic parameter choice.  The proof contains absolute constants, additive interval terms, and preprocessing costs that are absent from this proxy, so the table does not establish an input size at which the algorithm becomes faster.

\section{Conclusion}
\label{sec:conclusion-new}

The primitive-root trace problem is a local two-target problem.  Trace-admissibility is exactly equivalent to the pointwise bound
$|H_\gamma\cup H_{\gamma^2}|\le3\ell/4+1$, and the CM specialization
$\gamma=2$ supplies a GRH-conditional positive-density family.  The global
factorization argument developed here is independent of that conditional
trace-selection result.  It instead uses the residue-known
Harvey--Hittmeir coefficient lattice.

For the correct residue class of $p$ modulo $m$, the actual integer babystep
index in the Harvey--Hittmeir search has, modulo every selected prime
$\ell\mid m$, the form of a nondegenerate quadratic polynomial whenever
$\ell\nmid a$.  After a computable translation, every such quadratic image avoids the fixed set $F_\ell=\mathcal N_\ell\cap(1+\mathcal R_\ell)$, of size $(\ell-\chi(-1))/4$.  This quarter-density is optimal for a universal translated forbidden set that must handle both square classes of the quadratic leading coefficient.  By using only sufficiently large prime divisors of the auxiliary primorial, we can choose a nearby reduced-lattice vector that avoids all of them.  Gauss reduction then gives a regular case in which the quarter-sieve is used and a shorter unsieved case.  Arbitrary-point multipoint evaluation over $\Z/N\Z$ keeps the cost proportional to the number of retained exponents, and the unrestricted large-order theorem of Harvey and Hittmeir supplies the required element for the target order used here.

The endpoint choice $m=N^{1/10-o(1)}$ gives the deterministic bound in Theorem~\ref{thm:factorization-main}.  Its secondary factor improves the general Harvey--Hittmeir factorization bound~\cite{HarveyHittmeir2022} and, on balanced semiprimes, has a smaller secondary factor than the rank-$3$ bound with logarithmic exponent $13/5$~\cite{GaoFengHuPan2025}.  The factor $\sqrt{t\rho}$ in the search estimate comes from balancing the regular and skinny cases.  To improve the secondary constant within this method, one would need a better estimate for the skinny case or a different way to combine the two cases.

\end{document}